\documentclass{article}

\usepackage{graphicx}
\usepackage{hyperref}
\usepackage{amsmath}
\usepackage{amsfonts}
\usepackage{color}
\usepackage{fullpage}

\begin{document}

\title{A 2-Level Domain Decomposition Preconditioner for KKT Systems with Heat-Equation Constraints\footnote{To appear in the proceedings of   \emph{The 27th International Conference on Domain Decomposition Methods (DD27)}}
}
\date{}
\author{Eric C. Cyr\footnote{Eric C. Cyr: eccyr@sandia.gov}, Sandia National Laboratories, Albuquerqeue, NM}

\maketitle

\begin{abstract}
Solving optimization problems with transient PDE-constraints is computationally costly due  to the number of nonlinear iterations and the cost of solving large-scale KKT matrices. These
matrices scale with the size of the spatial discretization times the number of time steps. We propose a new
two level domain decomposition preconditioner to solve these linear systems when constrained by the heat equation. Our approach
leverages the observation that the Schur-complement is elliptic in time, and thus amenable to classical
domain decomposition methods. Further, the application of the preconditioner uses existing time integration routines to
facilitate implementation and maximize software reuse. The performance of the preconditioner is examined in an empirical study
demonstrating the approach is scalable with respect to the number of time steps and subdomains.
\end{abstract}

\section{Introduction}

This paper develops a new domain-decomposition method for solving the KKT system with heat-equation constraints. This effort is driven by the quadratic optimization problem of the form
\begin{equation} \label{Eq:quadratic-opt}
\begin{aligned}
\underset{z}{\min} \quad & \frac12 \int_0^T \|u-\tilde{u}\|_{L^2(\Omega)}^2\, dt +  \frac{\omega}{2} \int_0^T \|z\|_{L^2(\Omega)}^2 \, dt \\
\textrm{s.t.} \quad & \partial_t u - \nu\nabla\cdot\nabla u = z, & x\in\Omega \subset \mathbb{R}^2, t\in[0,T] \\
               & u(x,t) = 0, \;\; x\in{\partial\Omega}, t\in[0,T], &  u(x,0) = u_0(x), \;\; x\in\Omega 
\end{aligned}
\end{equation}
This quadratic PDE-constrained optimization problem finds a control $z$ such that the solution $u$ to the heat equation matches the target $\tilde{u}$.
The spatial domain is $\Omega$, the time interval is $[0,T]$, and the heat conductivity is $\nu$. Uniform homogenous boundary conditions are assumed for all time, and the initial condition is prescribed by $u_0$.

Many nonlinear methods use a series of quadratic approximations of the form represented by Eq.~\ref{Eq:quadratic-opt} to solve PDE-constrained optimization problems (see for instance sequential quadratic programming methods~\cite{byrd2008inexact,heinkenschloss2014matrix,ulbrich2007generalized}). There have been several studies focused on developing scalable preconditioners for the saddle-point system that arises from the first-order necessary conditions. Often preconditioners for saddle-point systems take the form of approximate factorization block preconditioners~\cite{benzi2005numerical}. These were explored for KKT systems in~\cite{benzi2011preconditioning,biros2005parallel}. Our work relies heavily on the block preconditioners from the Wathen group~\cite{pearson2012regularization,pearson2014preconditioners,rees2010optimal}. 

This effort focuses on transient PDE constraints where the size of the system scales with the number of spatial unknowns times the number of time steps, resulting in substantial computational effort.  To alleviate this, a number of efforts have proposed accelerating the time solve using adaptive space-time discretizations~\cite{langer2021unstructured,langer2021space}, parareal~\cite{maday2002parareal,ulbrich2007generalized,gander2020paraopt}, multigrid approaches~\cite{borzi2003multigrid,borzi2009multigrid,falgout2020trimgrit,gunther2019non,lin2022multilevel}, block preconditioning~\cite{pearson2012regularization,pearson2014preconditioners}, and domain decomposition methods~\cite{heinkenschloss2005linearquad}.  

Our approach is also built on block preconditioning ideas. A difference is that our technique exploits
an observation that the Schur-complement of the KKT system is elliptic in time (see~\cite{gander2016schwarz,lewis2005model}). This allows us to leverage existing two level domain decomposition approaches for elliptic systems to improve the parallel scalability of the block preconditioner. Good performance is achieved by algorithmic choices that ensure the forward and backward in time integrators can be applied on the fine level. 

\section{Discrete System and Block Preconditioner}\label{sec:discrete-system}

In this article the PDE in Eq.~\ref{Eq:quadratic-opt} will be discretized on a $2$D Cartesian grid using first order backward Euler in time, and a second order finite difference stencil in space. A row of the discrete
space-time system for the heat equation satisfies:
\begin{equation}\label{Eq:heat-fd}
u_{i,j}^{n+1} - u_{i,j}^n + \Delta t 
\nu
\left(
\frac{- u_{(i+1)j}^{n+1} + 2u_{ij}^{n+1}- u_{(i-1)j}^{n+1}}{\Delta x^{2}}
+\frac{- u_{i(j+1)}^{n+1} + 2u_{ij}^{n+1}- u_{i(j-1)}^{n+1}}{\Delta y^{2}}
\right)
=
 \Delta t  z_{ij}^{n+1}. 
\end{equation}
Here $i,j$ are the interior space indices defined over $1\dots n_x-1$ and $1\dots n_y-1$. The exterior indices are eliminated using the homogenous boundary conditions. The superscript time index $n$ runs from $0\ldots N_t$. Each $n$ is referred to below as a \emph{time-node}. The control variable $z$'s index matches the implicit index on $u$, therefore $z^{n+1}$ is associated with the $n$th time interval. For a single time interval, Eq.~\ref{Eq:heat-fd} rewritten in matrix form is
\begin{equation}\label{Eq:time-integration}
J_{(n+1)(n+1)} u^{n+1} + J_{(n+1) n} u_n + L_{(n+1)(n+1)} z^{n+1} = 0,
\end{equation}
and the global space-time system is
\begin{equation}
J u + L z = f.
\end{equation}
The right hand side $f$ includes contributions from the initial conditions.
The matrix $J $ is block lower triangular and
the matrix $L$ is block diagonal. 

The linear system whose solution solves the quadratic optimization problem from Eq.~\ref{Eq:quadratic-opt} is the celebrated KKT system 
$K \vec{u} = \vec{f}$ where
\begin{equation}\label{Eq:linsys}
K =
\begin{bmatrix}
M_u & & J^T \\
        & \omega M_z & L^T \\ 
J & L & 
\end{bmatrix}, \;
\vec{u}
=
\begin{bmatrix}
u \\ z \\ w
\end{bmatrix}, \;
\vec{f}
=
\begin{bmatrix}
f_u \\ f_z \\ f
\end{bmatrix}.
\end{equation} 
The final row is the discrete form of the PDE constraint, enforced by the Lagrange multiplier $w$.
We will also refer to $w$ as the adjoint solution. $M_u$ and $M_z$ are identity matrices scaled by $\Delta t \Delta x \Delta y$.
The matrix $K$ is a saddle point matrix, whose structure is frequently observed in numerical optimization. Many effective block preconditioners have been developed for this class of matrix~\cite{benzi2004preconditioner,benzi2005numerical,benzi2011preconditioning}. We focus on the block preconditioning approach developed by Wathen and collaborators for solving linearized PDE-constrained optimization problems~\cite{pearson2012regularization,pearson2014preconditioners,rees2010optimal}. 

We write a block LDU factorization of the matrix $K$
\begin{equation}
K =
\begin{bmatrix}
I & &  \\
        & I & \\ 
J M_u^{-1} & \omega^{-1} L M_z^{-1}& I 
\end{bmatrix}
\begin{bmatrix}
M_u & &  \\
        & \omega M_z &  \\ 
 &  & -S
\end{bmatrix}
\begin{bmatrix}
I & & M_u^{-1} J^T \\
        & I & \omega^{-1} M_z^{-1} L^T \\ 
 &  & I
\end{bmatrix}
\end{equation}
where the Schur-complement is
$S =  J M_u^{-1} J^T + \frac{1}{\omega} L M_z^{-1} L^T$. Following~\cite{pearson2012regularization}, $K$ is preconditioned using the block diagonal operator
\begin{equation}\label{Eq:block-prec}
P = 
\begin{bmatrix}
M_u & & \\
& \omega M_z & \\
& & \hat{S}
\end{bmatrix} \mbox{ where } \hat{S} = \hat{J} M_u^{-1} \hat{J}^T, \hat{J} = J + \omega^{-1/2} L.
\end{equation}
This preconditioner leverages the result in~\cite{murphy2000note}, and approximately inverts the block diagonal in the LDU factorization. 
 The matrix $\hat{J}$ used in the approximate Schur complement $\hat{S}$ is block lower triangular (similar to $J$), a fact that we will exploit below. The choice of $\hat{S}$ integrates the state Jacobian and the effects of the regularization parameter. In~\cite{pearson2012regularization, pearson2014preconditioners} and~\cite{schiela2014operator}, this approximation is developed and shown to lead to robust performance with respect to $\omega$.

\section{Two-Level Domain Decomposition Schur-Complement} \label{sec:domain-decomposition}

We propose a new domain decomposition approach for approximately inverting $\hat{S}$.
This is motivated by the observation that the operator $S$ is elliptic in time (see~\cite{lewis2005model} and~\cite{gander2016schwarz}). 
For simplicity, we show this discretely using only the term $J M_u^{-1} J^T$. Consider the 
ODE $\partial_t y = -y$ discretized over three time steps with forward Euler: $y^{n+1} - y^n + \Delta t y^n = 0$. With $M_u=I$, the Schur-complement $\hat{S}$ is
\begin{multline}
\begin{bmatrix}
1 \\
-1+ \Delta t & 1  & & \\
& -1+ \Delta t & 1  & \\
\end{bmatrix}
\begin{bmatrix}
1 & -1 + \Delta t  & \\
& 1 & -1 + \Delta t  \\
 &  &                    1 \\
\end{bmatrix} \\ = 
\begin{bmatrix}
1 & -(1 - \Delta t) & \\
-(1 - \Delta t)   & 2 (1-\Delta t) + \Delta t^2 & -(1 - \Delta t ) \\
& -(1 - \Delta t)   & 2 (1-\Delta t) + \Delta t^2 
\end{bmatrix}.
\end{multline}
Examining the second row it is clear the operator has a $1$D Laplacian stencil in time, with a positive perturbation on the diagonal. To take advantage of this ellipticity, we will apply existing domain decomposition approaches to the $\hat{S}$ operator. This \emph{ellipticity principle} enables scalable performance of a preconditioned Krylov method. We also impose an \emph{efficiency constraint} that the computational kernels in our preconditioner use the time integration method for the state and adjoint unknowns. 

The \emph{ellipticity principle} is realized by considering a restricted additive Schwarz (RAS) method with $N_D$ subdomains (see~\cite{Dolean2015}). Each subdomain contains the spatial unknowns associated with a subset of time steps. For instance, if there are $N_t = 9$ time steps, then for $N_D=3$ the subdomains contain time-nodes $\{1,2,3\}$,
$\{3,4,5,6\}$, and $\{6,7,8,9\}$ (the $0^\text{th}$ time-node is the excluded initial condition). Notice that the time-nodes are overlapped but the time steps are not. With these subsets, boolean operators $R_s$ are defined that restrict a space-time vector to the time-nodes in a subdomain, giving the RAS preconditioner
\begin{equation} \label{Eq:sas}
\hat{S}_{\text{RAS}}^{-1} = \sum_{s=1}^{N_D} R_s^T D_s \left(R_s \hat{J} M_{u}^{-1} \hat{J}^T R_s^T \right)^{-1} R_s
\end{equation}
where $D_s$ is the (boolean) partition of unity matrix (Defn.~1.11 of~\cite{Dolean2015}). 
RAS is known to lead to effective preconditioners for elliptic problems and can be extended to a multi-level schemes. While the \emph{ellipticity principle} is exploited in $\hat{S}_{\text{RAS}}^{-1} $, the explicit formation of the product $\hat{J} M_u^{-1} \hat{J}^T$ does not satisfy the \emph{efficiency constraint}.

To satisfy the \emph{efficiency constraint} note that the range of $\hat{J}^T R_s^T$ is nonzero on time-nodes in the $s^\text{th}$ subdomain and one time-node earlier. For instance, if the subdomain contains nodes $\{3,4,5,6\}$ then the range is nonzero on $\{2,3,4,5,6\}$. Let $Q_s$ be a new extended restriction operator whose action produces
a space-time vector for time-nodes that are nonzero in the range of $\hat{J}^T R_s^T$. Further, choose $Q_s$ so that
\begin{equation}\label{Eq:def-qs} 
Q_s = 
\begin{bmatrix}
W_s \\
R_s
\end{bmatrix}
\mbox{ and } Q_s Q_s^T = I.
\end{equation}
The operator $W_s$ restricts the space-time vector to the time-nodes contained in the earlier time step relative to the $s^{\text{th}}$ subdomain.
Because $Q_s$ is a restriction operator that represents the nonzero range of $\hat{J}^T R_s^T$, we have that $Q_s^T Q_s \hat{J}^T R_s^T = \hat{J}^T R_s^T$.
Recalling that $M_u$ is diagonal, we can rewrite the subdomain solve in $\hat{S}_{RAS}^{-1}$ as
 \begin{equation} \label{Eq:subdomain-equality}
R_s \hat{J} M_u^{-1} \hat{J}^T R_s^T= R_s \hat{J} Q_s^T (Q_s M_u^{-1} Q_s^T) Q_s \hat{J}^T R_s^T
\end{equation} 
Using the constraint in Eq.~\ref{Eq:def-qs}, we have the additional identities $R_s =  R_s Q_s^T Q_s$ and $R_s^T =  Q_s^T Q_s R_s^T$.
This permits a final rewrite of the operator in Eq.~\ref{Eq:subdomain-equality} 
\begin{equation}
R_s \hat{J} M_u^{-1} \hat{J}^T R_s^T        = R_s Q_s^T \hat{J}_s M_s^{-1} \hat{J}_s^T Q_s R_s^T
\end{equation}
where $\hat{J}_s = Q_s \hat{J} Q_s^T$ and $M_s^{-1} = Q_s M_u^{-1} Q_s^T$. The inverse action of $\hat{J}_s M_s^{-1} \hat{J}_s^T$ is easily computed
in a matrix free way on the time-nodes in the extended subdomain. Motivated by this equivalence, we define a new one-level preconditioner
\begin{equation} \label{Eq:sasq}
\hat{S}_{\text{RASQ}}^{-1} = \sum_{s=1}^{N_D} R_s^T D_s \left(R_s Q_s^T \hat{J}_s^{-T} M_s \hat{J}_s^{-1} Q_s R_s^T \right)  R_s.
\end{equation} 
The term in parentheses is different from the term inverted in Eq.~\ref{Eq:sas}. The difference is that the inverse computed in Eq.~\ref{Eq:sas} is constrained to have a zero initial condition outside of the subdomain. This revised operator satisfies our \emph{efficiency constraint} as computing $\hat{J}_s^{-1}$ and $\hat{J}_s^{-T}$ is done using the time integration method.

To obtain scalability with respect to the number of subdomains a coarse grid correction is required. Again leveraging the elliptic
nature of the Schur-complement we consider the Nicolaides coarse space developed for solving the Poisson problem~\cite{Nicolaides1987}. 
Following the presentation in~\cite{Dolean2015}, the Nicolades coarse space is defined as 
\begin{equation}
Z = \begin{bmatrix}
D_1 R_1 \Phi_1 &   & \ldots & \\
& D_2 R_2 \Phi_2 & \ddots  & \\
\vdots & \ddots & \ddots  & \\
          & \ldots &   & D_{N_D} R_{N_D} \Phi_{N_D}
\end{bmatrix}
\mbox{ where }
\Phi_s = 
\begin{bmatrix}
1 &  -1        \\
1 & -1+\frac{2}{N_s}\\
\vdots & \vdots \\
1  & -1 +2\frac{N_s}{N_s}
\end{bmatrix}
\end{equation}
The columns of $\Phi_s$ form a constant and linear basis over the $N_s$ subdomain time-nodes.
The coarse restriction $R_0 = Z^T$ is used in the definition of the coarse operator
$\hat{S}_0 = R_0 \hat{J} M_u^{-1} \hat{J}^T R_0^T$.
The coarse solve is applied in a multiplicative way
\begin{equation}
\hat{S}_{2-\text{level}}^{-1} = \hat{S}_{\text{RASQ}}^{-1} R_0^T \hat{S}_0^{-1} R_0.
\end{equation}
Due to the structure of $R_0$ the coarse operator $\hat{S}_0$ can be constructed in parallel. This does represent a violation of the \emph{efficiency constraint} to be addressed by future work.

\section{Numerical Experiments} \label{sec:num-experiments}

To demonstrate this approach we discretize the quadratic optimization problem from Eq.~\ref{Eq:quadratic-opt} as described in Sec.~\ref{sec:discrete-system}. 
The $2$D spatial domain is $\Omega = (0,1) \times (0,2) \subset \mathbb{R}^{2}$, and the time domain is $[0,1]$. The initial conditions and target solutions are
\begin{equation}
u_0(x,y) = -xy(x-1)(y-2), \quad \tilde{u}(x,y,t) = \sin(2.0\pi t)\sin(2.0\pi x)\sin(2.0\pi y).
\end{equation}
The regularization parameter $\omega$ varies over five orders of magnitude.
Experiments were
run with $9\times 9$, $17\times 17$, and $33\times 33$
mesh points. Qualitatively, variation with number of spatial points was not a factor in the convergence. This is not surprising as the implicit operator in space is inverted with a direct solve. 
As a result, the computations below are all for the case of $17\times 17$ mesh points. Recall that homogeneous boundary conditions
are removed, giving $15\times 15$ unknowns in each time step. The linear system $K \vec{u} = \vec{f}$ (Eq.~\ref{Eq:linsys}) is solved using right preconditioned GMRES from PyAMG~\cite{PyAMG2022} iterated until a relative residual tolerance of $10^{-6}$ is achieved.

\begin{figure}
\centering{
\includegraphics[width=0.32\textwidth]{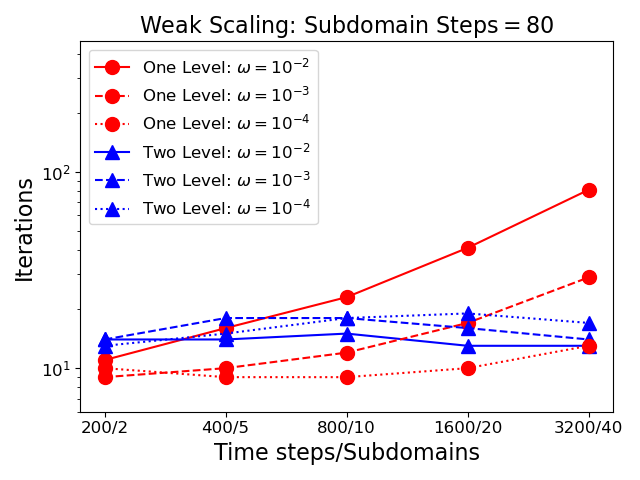}
\includegraphics[width=0.32\textwidth]{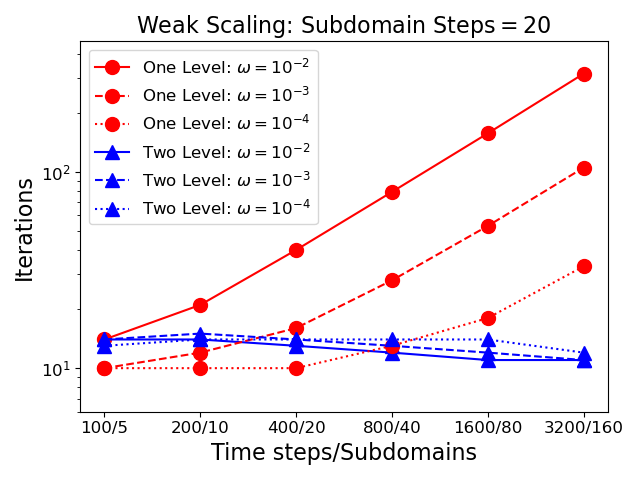}
\includegraphics[width=0.32\textwidth]{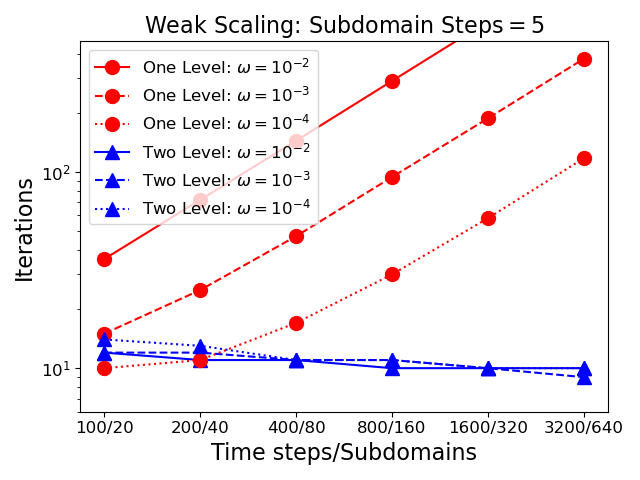}

}
\caption{Three weak scaling studies for different numbers of time-steps per subdomain. The two level scheme (triangles) has flat iteration
counts for regardless of the number of time steps, the subdomain size, and the regularization parameter. Asymptotically
the one level method shows a strong dependence with respect to the number of subdomains and time steps.}
\label{fig:weak}
\end{figure}

\begin{figure}
\centering{
\includegraphics[width=0.32\textwidth]{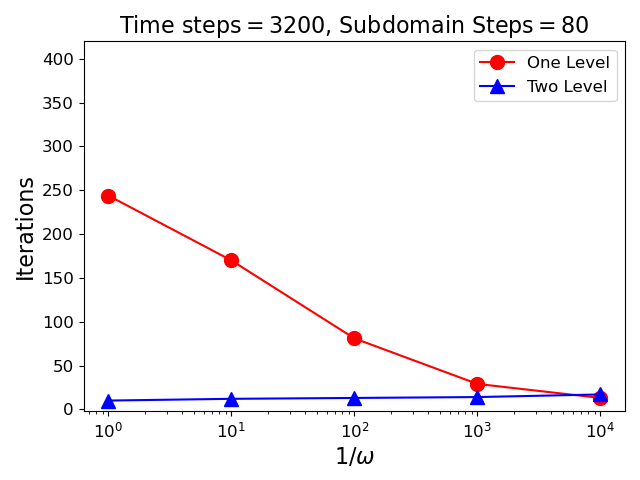}
\includegraphics[width=0.32\textwidth]{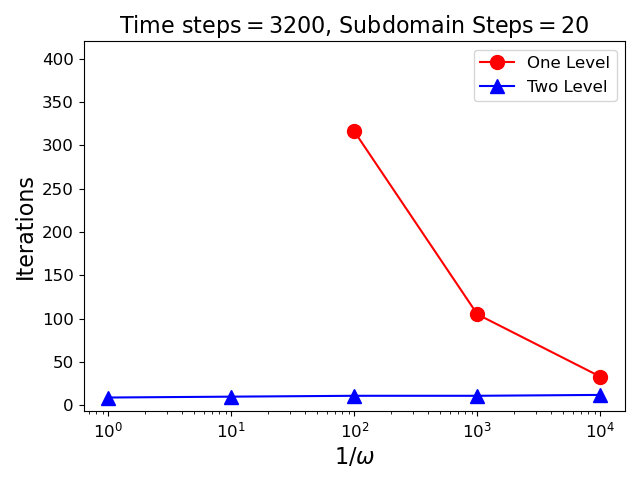}
\includegraphics[width=0.32\textwidth]{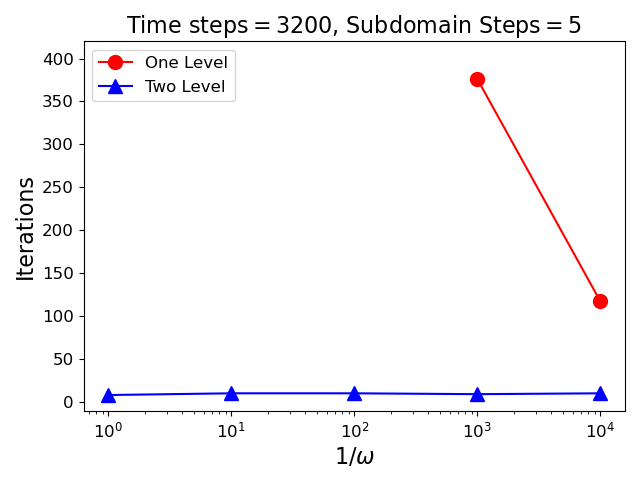}
}
\caption{This plot demonstrates the robustness of the two level scheme (triangle markers) with respect to the regularization parameter $\omega$. 
Note that for many cases the one level scheme (circle markers) did not converge in the $420$ iterations (the maximum allowed), thus those
values are omitted.}
\label{fig:omega}
\end{figure}

Figure~\ref{fig:weak} presents three weak-scaling studies ranging from $100$ to $3200$ time steps. The number of time steps per subdomain is fixed at $80$ in the left plot,
$20$ in the center plot, and $5$ in the right plot. For the case of $80$ steps, the fewest number of time steps is $200$ (the minimum number of time steps evenly divisible by $80$ in the chosen sequence).
The plots, show the number of iterations as a function of time step count for GMRES preconditioned using $P$ from Eqn~\ref{Eq:block-prec} with Schur complement approximations $\hat{S}_{\text{RASQ}}^{-1}$ for the one level case (circle markers), and $\hat{S}_{2-\text{level}}^{-1}$ for the two level case (triangle markers). Different values for the regularization parameter $\omega$ are indicated using solid ($10^{-2}$), dashed ($10^{-3}$) or dotted ($10^{-4}$) lines.  These plots demonstrate that the performance of the two level method is independent of both
the number of subdomains, and the number of time steps. Further, independence holds regardless of the value of the regularization parameter. As anticipated the one level method has substantial growth with the number of time steps, and variability with the regularization parameter. However, it is worth noting that dependent on the number of subdomains and the size of the regularization parameter the one level method may be faster despite its lack of scalability. For instance, when using $40$ subdomains and a regularization parameter of $10^{-4}$ the one-level method takes the same number of iterations but lacks the synchronization and added cost of the two level method.

The scaling with respect to the regularization parameter is investigated in Figure~\ref{fig:omega}. In these plots the preconditioned iteration counts are plotted as a function of the inverse regularization parameter. Data points are excluded when the number of iterations exceeded the maximum iteration count for GMRES (in this case $420$). 
Here again the two level method scales well, yielding essentially flat
iteration counts as a function of the regularization parameter. The one level method shows strong dependence on $\omega$, though it improves dramatically for smaller values. 

\section{Conclusion} \label{sec:conclusion}
In this paper, motivated by results in block preconditioning and the elliptic-in-time nature of the KKT system, we develop a two level domain decomposition preconditioner that facilitates a parallel-in-time solver for the discrete optimality system constrained by the heat equation. While limited in their breadth, initial results for this approach show excellent scalability with respect to the number of time steps, subdomains, and the regularization parameter. Future work will focus on achieving improved scaling by including more levels in the hierarchy, and applying this technique to a broader class of problems and discretizations.

\section*{Acknowledgements}
The author is indebted anonymous reviewers whose comments resulted in a significant strengthening of the paper.
The author also acknowledges support from the SEA-CROGS project and the Early Career program, both funded by the DOE Office of Science.
This article has been authored by an employee of National Technology \& Engineering Solutions of Sandia, LLC under Contract No. DE-NA0003525 with the U.S. Department of Energy (DOE). The employee owns all right, title and interest in and to the article and is solely responsible for its contents. The United States Government retains and the publisher, by accepting the article for publication, acknowledges that the United States Government retains a non-exclusive, paid-up, irrevocable, world-wide license to publish or reproduce the published form of this article or allow others to do so, for United States Government purposes. The DOE will provide public access to these results of federally sponsored research in accordance with the DOE Public Access Plan https://www.energy.gov/downloads/doe-public-access-plan. This paper describes objective technical results and analysis. Any subjective views or opinions that might be expressed in the paper do not necessarily represent the views of the U.S. Department of Energy or the United States Government. Sandia report number: SAND2023-03175O.


\bibliographystyle{plain}
\bibliography{main}

\end{document}